\documentclass[a4paper,10pt]{article}

\usepackage[a4paper, textwidth=160mm, marginratio=1:1]{geometry}

\usepackage{amsmath, amssymb, amscd, amsthm}

\usepackage{amsfonts}

\usepackage[all]{xy}

\usepackage{graphics}

\usepackage{graphicx}

\usepackage{epsfig}

\usepackage{latexsym}

\theoremstyle{plain}

\newtheorem{thm}{Theorem}[section]

\newtheorem{theorem}[thm]{Theorem}

\newtheorem{prop}[thm]{Proposition}

\newtheorem{proposition}[thm]{Proposition}

\newtheorem{coro}[thm]{Corollary}

\newcommand{\lra}{\longrightarrow}

\newcommand{\RR}{\mathbb{R}}

\newcommand{\NN}{\mathbb{N}}

\newcommand{\ZZ}{\mathbb{Z}}

\newcommand{\OO}{\mathcal{O}}

\newcommand{\K}{\mathcal{K}}

\def\bbR{\mathbb{R}}

\def\bbC{\mathbb{C}}

\def\bbZ{\mathbb{Z}}

\def\bbC{\mathbb{C}}

\def\z2{\mathbb{Z}/2}

\def\ztwo {\mathbb{Z}/2\mathbb{Z}}

\newcommand{\sphere}[1]{\mathbb{S}^{#1}}

\def\cliff#1{\mathcal{C}{\ell}_{#1}}                

\def\ccliff#1{\mathbb{C}{\ell}_{#1}}

\newcommand{\hilb}{\mathcal{H}}

\newcommand{\chilb}{\mathcal{L}}                    

\newcommand{\hilbR}[2]{\mathcal{R}^{#2}_{#1}}

\newcommand{\hhilbR}[2]{{\mathfrak{R}}^{#2}_{#1}}

\newcommand{\conjugate}[1]{\overline{#1}}

\def\KO {\mathrm{KO}}

\def\KR {\mathrm{KR}}

\def\KU {\mathrm{KU}}

\def\kr {\mathrm{kr}}

\def\ko {\mathrm{ko}}

\def\K {\mathrm{K}}

\def\krpqip1{\kr_{p,q}(i+1)}

\newcommand{\HFF}[2]{\hat{\mathbb{F}}^{#2}_{#1}}

\newcommand{\innfR}[2]{\mathcal{R}Inf_{#1}^{#2}}

\newcommand{\hinnfR}[2]{\mathcal{R}I\widehat{n}f_{#1}^{#2}}

\newcommand{\hinf}{I\widehat{n}f}

\newcommand{\flag}{\mathcal{F}lag}

\newcommand{\nerve}{\mathfrak{NC}}

\def\colim {\mathrm{colim}}

\newcommand{\Rinfty}{\overline{\mathbb{R}}}

\def\pminusone#1{p^{-1}(#1)}

\newcommand{\proofbox}{\begin{flushright}$\Box$\end{flushright}}

\def\domperp#1{\mathcal{D}(#1)^{\perp}}

\newcommand{\innf}[2]{Inf_{#1}^{#2}}

\newcommand{\hinnf}[2]{I\widehat{n}f_{#1}^{#2}}

\newcommand{\Esp}[2]{Esp(#1;#2)}

\def\gradedtensor{\hat{\otimes}}

\def\dom#1{\mathcal{D}(#1)}

\newcommand{\D}{\mathcal{D}}

\newcommand{\homeo}{\approx}

\newcommand{\fingrass}[2]{\mathcal{G}r_{#2}(#1)}

\def\set#1{\{#1\}}                            

\def\tensor{\otimes}                               

\newcommand{\Conflinfin}[2]{{\it C}\widehat{{\it on}}{\it f}_{#1}^{#2}}

\def\HEE#1#2{\widehat{\mathbb{E}}_{#1}^{#2}}

\newcommand{\bottzeroone}{\beta_{0,1}}

\newcommand{\bottonezero}{\beta_{1,0}}

\newcommand{\bottoneone}{\beta_{1,1}}

\def\covermap {\mathrm{\varphi}}

\def\sm {\smallskip}

\title{An Operator model for connective $K$-theory with reality}

\author{G. Gaudens, E. Markert}

\date{2011}

\begin{document}

\maketitle






\begin{abstract}

We construct an explicit model for the connective cover of the spectrum of $K$-theory with reality.

This model is a $\ztwo$-equivariant commutative symmetric ring spectrum which

identifies with another one due to D. Dugger.

\end{abstract}

\section{Models for $\K$-theory with reality}

\subsection{Periodic $\K$-theory with reality}

$\K$-theory is a generalized cohomology theory that comes in
various flavours; the most famous ones are certainly the real
$\K$-theory $\KO$ and the complex $\K$-theory $\KU$. Many explicit
models for a spectrum representing these cohomology theories are
known (\cite{atiyahsinger, HST}). Both $\KU$ and $\KO$ are
contained in a $\ztwo$-equivariant cohomology
theory called $\K$-theory with reality (denoted by $\KR$)
discovered by Atiyah \cite{atiyah}. For instance, Bott periodicity
for $\KR$ (the 1-1 periodicity theorem) implies at once
by a clever argument  both real and complex Bott periodicity
theorems.

\sm

The following three papers provide the basis for this exposition:

\begin{itemize}

\item[-]  \cite[Theorem 3]{HST}, where a nice operator model for the periodic real $\K$-theory spectrum is defined, and models for the spectra of the complex $\K$- theory and of $\K$-theory with reality are suggested,

\item[-]  \cite{markert} where such a model for the real \emph{connective} real $\K$-theory spectrum $\ko$ is given,

\item[-]  \cite{dugger} where a model for connective $\K$-theory with reality is constructed by an abstract localization method.

\end{itemize}

Here, a \emph{nice} operator model for $\KO$ refers to a model that forms a \emph{commutative symmetric ring spectrum}, whose product actually  lifts the natural product on $\K$-groups. This was first observed by M. Joachim \cite{joachimKO, joachimKU}, while another approach appears in \cite{mitchener}.

\sm

Our aim in this work is to provide an operator model for the
spectrum $\KR$ as outlined in \cite{HST} and for the corresponding
connective version  $\kr$ following the method in \cite{markert}.
We show that this model for $\kr$ coincides with the one
constructed in \cite{dugger}.
Finally, we know from \cite{BR1, BR2} that highly structured ring structures for $\KU$, $\KO$ and their connective covers exist and are essentially unique. We give here an explicit model for such a cover in the case of $\KR$.  So our model for $\KR$ provides also in particular explicit such coverings for $\KO$ and $\KU$. In view of these strong uniqueness properties and the universality of Dugger's construction, it is somehow  surprising to find such an explicit model with the required properties, {\it i.e.} it is a connective cover as a \emph{ring spectrum} \emph{and} it has the `right' equivariant homotopy type.

\sm

This altogether suggests that the operator setting emphasized in \cite{HST, markert} is particularly suited to provide explicit models for $\K$-theory `in the large', for instance in the case of  $G$-equivariant connective $K$-theory \cite{greenlees}. We consider the case of $\kr$ as a non trivial test case for constructing  general connective models for topological $\K$-theory spaces, which will be the subject of future work by the authors.

\sm {\it Acknowledgments.~} The first author wishes to thank the University of Bonn and the Max Planck Institut of Bonn for
their support in elaborating this work.

\subsection{Conventions and notations}

We will use the notation (possibly decorated by a subscript)
$\hilb$ to denote a separable, infinite-dimensional $\ztwo$-graded
real Hilbert space, and use the notation $\chilb$ for a complex
version.

\sm

For positive integers $p$, $q$, we denote by $\cliff{p,q}$ the
Clifford-algebra on $\bbR^{p+q}=\bbR^p\oplus\bbR^q$ with the
positive definite standard quadratic form on $\bbR^p$ and the
negative definite standard form on $\bbR^q$. That is,
$\cliff{p,q}$ is the quotient of the free associative
$\bbR$-algebra on generators $e_i$, $f_j$, $i=1, \ldots, p$,
$j=1,\ldots, q$ by the two-sided ideal generated by the relations
$e_i^2=-1$, $f_j^2=1$ and the relations given by the fact that any
two different generators anti-commute. We shorten notation by
using $\cliff{n}:=\cliff{n,0}$ and $\cliff{-n}:= \cliff{0,n}$ for
positive integers $n$. The complex Clifford algebras
$\ccliff{p,q}$ are defined analogously, on $\bbC^p\oplus\bbC^q$.
The algebras $\ccliff{p,q}$ and $\cliff{p,q}\tensor \bbC$ are
canonically isomorphic.

We use the same short notations in the complex case.

\sm

The Clifford algebras are $\ztwo$-graded with generators in odd
degree, the grading involution will be denoted by $\epsilon$. Note
that there is a correspondence between $\ztwo$-graded Clifford
structures and ungraded ones via an index shift by one (the
even part of $\cliff{n}$ is isomorphic to $\cliff{n-1}$). In this
paper we work in the graded setting.

\sm

Furthermore, there is a grading-preserving involution $c$ on
$\cliff{p,q}$, sending the $(f_j)_{ij=0,\ldots, q}$ to their
negatives. This involution extends uniquely to an antilinear
involution $\theta$ on $\ccliff{p,q}$. This turns the algebras
$\cliff{p,q}$ and $\ccliff{p,q}$ into \emph{Real modules} over
themselves. A Real module $M$ for $\cliff{p, q}$ is a complex
representation endowed with an antilinear involution $c$, such
that the Clifford action of $x$ on $m\in M$ satisfies
$c(x.m)=\theta(x).c(m)$. One can extend a Real module for
$\cliff{p, q}$ uniquely to a Real module for $\ccliff{p, q}$. All
three notions of $\cliff{p, q}$-module, Real module for $\cliff{p,
q}$ and Real module for $\ccliff{p, q}$ produce isomorphic module
categories. For our $K$-theory models, we will use the notion of
Real module for $\ccliff{p, q}$.

We refer the reader to \cite{lawsonmichelson} for the general theory of Clifford algebras and their modules.

However one aspect we want to emphasize in advance is the
isomorphism of module categories between left $\cliff{n}$-modules
and right $\cliff{-n}$-modules for all $n\in\bbZ$.

We define
$$\hilbR{p, q}{}:= (\hilb\gradedtensor \ccliff{1})^{\gradedtensor
_{\bbC}p}\gradedtensor _{\bbC}(\hilb\gradedtensor
\ccliff{-1})^{\gradedtensor _{\bbC}q}\cong \chilb_p \gradedtensor
_{\bbC} \chilb_{-q}\cong \chilb \gradedtensor _{\bbC}
\chilb\gradedtensor_{\bbC}\ccliff{p, q}$$ where
$\gradedtensor$ denotes the graded tensor product of graded
Hilbert spaces (also, suitably completed).

The involution $c$ on $\ccliff{p,q}$ extends uniquely to an
antilinear involution $\theta$ on $\hilbR{p, q}{}$. This
involution turns $\hilb_{p, q}$ into a Real module for
$\ccliff{p,q}$.

\sm All future mentions about an action of $\ztwo$ in the
following will be correspond to the action induced by the
involution $\theta$.

\sm

The generators of $\ccliff{p,q}$ act on   $\hilb_{p,q}$ by
orthogonal transformations. Moreover   $\hilb_{p,q}$ is a
\emph{universe} for the  Clifford algebra $\ccliff{p,q}$ in the
sense that it contains each irreducible Real module infinitely
many times. Note that any two such universes are non-canonically
isomorphic. This does not matter in the homotopy category, because
the automorphisms of such a universe form a contractible space, by
the natural generalization of a theorem of Kuiper \cite{kuiper}.

\sm

We will need to work with a Real version of these Clifford
universes and we will also extend the concept to bi-modules. For
instance the Real universe $\hilbR{p+k,q+l}{}$ for
$\ccliff{p+k,q+l}{}$ is also a Real universe for the smaller
Clifford algebras $\ccliff{p,q}$. It will be important to keep
track of two different Clifford actions by separating them, that
is by letting the algebra $\ccliff{p,q}{}$ act on
$\hilbR{p+k,q+l}{}$ from the right as usual, and thinking of the
action of $\ccliff{k,l}{}$ as an action of the opposite algebra
$\ccliff{l,k}$ {\it from the left}. That is, the Hilbert universe
$\hilbR{p+k,q+l}{}$ becomes a
$\ccliff{l,k}-\ccliff{p,q}$-bimodule, which is simultaneously a
right $\ccliff{p,q}$-universe and a left $\ccliff{l, k}$-universe.
To express this in the notation we write $\hilbR{p,q}{k,l}$.

In the following we will consider
$$\hhilbR{p,q}{}:=\hilbR{p,q}{q,p}.$$
This is a universe for $\ccliff{q,p}-\ccliff{p,q}$ Real
bi-modules.

By convention, we set $\hhilbR{0,0}{}:=\chilb \gradedtensor _{\mathbb{C}}\chilb$.

\subsection{Spaces for $K$-theories}

Using the notations of \cite{HST} and \cite{markert}, we say that
an operator $G$ on a graded Hilbert space is of type $\hinf$ if it
is self-adjoint on a closed domain $\dom{G}$ and has finite rank
resolvent. Note that self-adjoint in general means densely
defined; however the operators in $\hinf$ are self-adjoint
operators on a finite-dimensional closed domain. If one requires
the resolvent to be only compact one obtains a larger class
$\innf{}{}$ of operators which are self-adjoint on a possibly
infinite-dimensional closed domain (where they are densely
defined). The spectral theorem implies that both types of
operators have a spectral decomposition into pairwise orthogonal,
finite-dimensional eigenspaces with real eigenvalues which can -
in case of $\innf{}{}$ - accumulate at infinity (the base point in the one-point compactification
$\Rinfty$). Even though the domain is part of the information, we
will in general suppress notation of the domains. The topology on
the spaces of operators of type $\innf{}{}$ is the coarsest
topology making functional calculus with functions in $C^c(\RR)$
continuous; here $C^c(\RR)$ denotes functions with compact
support. This topology allows us to use the spectral decomposition
in an intuitive way; we can think of the operators as
configurations of discrete points (eigenvalues) on the real line
with coefficients (eigenspaces) given by pairwise orthogonal,
finite-dimensional subspaces of the  background Hilbert space.

Now configuration points can be moved continuously (by functional calculus), and coefficients can be changed continuously within the
space of projection operators onto modules as above. Whenever two
configuration points collide, their coefficients add in algebraic
sum, whenever points move to $\infty$, they disappear from the
configuration. This is explained in detail in \cite{HST, markert},
including precise definitions of the configuration spaces and
their topology. Under this topology the spaces $\innf{}{}$ and
$\hinf$ are homotopy equivalent via the natural inclusion.
We will work throughout with the finite-rank version $\hinf$.
The space of $\hinf$-operators is naturally pointed, with base
point the zero domain operator.

\sm

Let $\hinnfR{p,q}{}$ denote the space of \emph{ odd, }
$\ccliff{p,q}$\emph{-linear} operators on $\hhilbR{p,q}{}$ which
are of type $\hinf$. This is a $\ztwo$-space under the action
induced by the involution $\theta$, sending an operator $G$ to its
conjugate $\theta G \theta$. This operator has the same
eigenvalues as $G$; but where $G$ has eigenspace $G_\lambda$ at
$\lambda$, the conjugate has eigenspace $\theta(G_\lambda)$.
We set
$$\KR_{p, q}(\hilb):= \hinnfR{p,q}{}.$$
Note that here we only use the right Clifford-action on the
universe $\hhilbR{p,q}{}=\hilbR{p,q}{q,p}$.

The chosen Hilbert universe $\hilb$ fixed in the beginning is suppressed in the notation, since isomorphisms of Hilbert universes induce homeomorphisms on the operator spaces.
The following is implicit in \cite{HST}:
\begin{prop}
For $p, q\geq 0$ and a $\ztwo$-CW-complex $X$, there is a natural
bijection of pointed sets
$$\KR_{p, q} (X)\cong [X, \KR_{p, q}]$$
where homotopy classes are through equivariant
homotopies.
\end{prop}

Notice that our conventions differ slightly from those of
\cite{HST}.

Firstly, we use right Clifford actions instead of the
corresponding left actions of the opposite algebras. Furthermore,
our spaces $\KR$ are defined as the finite-rank resolvent version
$\hinf$ of the operator spaces, whereas in \cite{HST} the authors
use the compact resolvent version $\innf{}{}$ and show that this
is homotopy equivalent to the latter (by an {\it ad hoc}
adaptation of \cite[Proposition 29]{HST}).

\sm

We now define spaces of operators $\hinnfR{p,q}{k,l}$ which will be shown to be connective covers of $\hinnfR{p,q}{}$, in the sense that there are equivariant homotopy equivalences
$$\hinnfR{p,q}{}\stackrel{\sim}{\longrightarrow} \Omega^{k,l}\hinnfR{p,q}{k,l}$$
where the spaces $\hinnfR{p,q}{k,l}$ are highly connected (Proposition \ref{overview}).

We set:
$$\hinnfR{p,q}{k,l}:=\set{G\in\hinnfR{p,q}{}(\hilbR{p,q}{k,l}) | G^2 \text{ is } \ccliff{l,k}-\ccliff{p,q}\text{-linear}}.$$

\sm

We define our spaces for connective $\K$-theories as
$$\kr_{p, q}(\hilb):=\hinnfR{p,q}{q,p}.$$
As before, we will usually suppress $\hilb$ from the notation.
Observe however that the definition still makes sense if we
replace $\hilb$ by a finite dimensional Hilbert space.

This is again a space of operators on the universe $\hhilbR{p,q}{}$ for all $p,q\geq 0$, this time the operators are defined using both Clifford actions.

\sm

There is an obvious map
$$\covermap_{p, q} : \kr_{p, q}   \to \KR _{p,q}$$
obtained by observing that $\hhilbR{p, q}{}$ besides being a
universe for $\ccliff{p, q}-\ccliff{p,q}$-bimodules, is a also a
universe for right $\ccliff{p, q}$-modules, hence $\kr_{p, q}$ is
a \emph{subspace} of $\KR_{p, q}$ for all $p, q\geq 0$.

So the spaces of operators on both source and target of
$\covermap$ do not live on the same Hilbert space. We leave it
to the reader to keep track of this subtlety, as it does not alter
the statements.

\subsection{Spectra for $\K$-theories}

\label{def:spaces:K}

\subsubsection{Symmetric structure and multiplication}

$\hilbR{p, q}{k,l}$ has an action of the group $\Sigma_p \times
\Sigma_q \times \Sigma_k \times \Sigma_l$ by permuting the tensor
factors. The symmetric group $\Sigma_p \times \Sigma_q$ acts on
$\hhilbR{p, q}{}= \hilbR{p, q}{q,p}$ diagonally.

\sm If $\varphi$ is in $\hinnfR{p,q}{k,l}$  and $ \psi$ is in  $\hinnfR{p',q'}{k',l'}$, then the operator
$$\varphi\star \psi:=\varphi  \gradedtensor {\mathrm I} + {\mathrm I} \gradedtensor \psi$$
is in $\innfR{p+p',q+q'}{k+k', l+l'}$.  This follows from
spectral analysis.

\sm Hence we get for all $p,p', q,q'\geq 0$, a commutative diagram of $\Sigma_p \times \Sigma_q \times \Sigma_{p'} \times \Sigma_{q'}$-equivariant maps

\begin{eqnarray}
\label{com:mult:cover}
\xymatrix{
\kr_{p,q} \wedge \kr_{p', q'} \ar[d]^{\covermap_{KR}} \ar[r]^{\star}&\kr_{p+p', q+q'}\ar[d]^{\covermap_{KR}}\\
\KR_{p,q}  \wedge \KR_{p', q'}  \ar[r]^{\star} &\KR_{p+p', q+q'}&.
}
\end{eqnarray}

The maps $\star$ are $\ztwo$-equivariant in the sense that the diagrams
\begin{eqnarray}
\label{equiv:mult} \xymatrix{
\KR _{p,q} \wedge \KR_{p',q'} \ar[r]\ar[d]^{\theta\wedge \theta}&  \KR_{p+p',q+q'}\ar[d]^{\theta}&\kr _{p,q} \wedge \kr_{p',q'} \ar[r]\ar[d]^{\theta\wedge \theta}&  \kr_{p+p',q+q'}\ar[d]^{\theta}\\
\KR_{p,q} \wedge \KR_{p', q'} \ar[r]& \KR_{p+p',q+q'}&\kr_{p,q}
\wedge \kr_{p', q'} \ar[r]& \kr_{p+p',q+q'} }
\end{eqnarray}
are commutative diagrams of $\Sigma_p \times \Sigma_q \times
\Sigma_{p'} \times \Sigma_{q'}$-spaces.

\subsubsection{Unit maps}\label{unitmaps}

We will use the definition of an equivariant symmetric ring
spectrum by \cite{mandell}.

\sm

Let $\mathbb{S}^{\bullet, \bullet}$ be the $\ztwo$-equivariant
commutative symmetric ring spectrum of spheres, where
$\mathbb{S}^{p,q}$ is the one-point compactification of the
$\ztwo$-representation $\mathbb{R}^{p,q}$ which is $p$ times the
trivial representation plus $q$ times the sign representation. We
define  $\ztwo$-equivariant maps
$$
\eta_{\kr}:\mathbb{S}^{\bullet, \bullet} \lra \kr_{\bullet,
\bullet}, ~~~\eta_{\KR}:\mathbb{S}^{\bullet, \bullet} \lra
\KR_{\bullet,\bullet}~~.
$$

\sm

We proceed as follows (compare \cite{joachimKO, HST}). Let $v$ be an
element in $\RR ^{p,q}$. We can see $v$ as an element of the
Clifford algebra $\ccliff{p,q}$, and let it act via  Clifford
multiplication {\it from the left} by $v$  on the Clifford
universe $\hilbR{p, q}{q,p}$. We denote this operator by $L_v$.
Recall here that $\hilbR{p,q}{q,p}$ is a
$\ccliff{p,q}-\ccliff{p,q}$-bimodule (a $\ccliff{p+q,p+q}$-right
module, respectively). We are considering $v$ as elements of the
algebra acting from the left (or as the corresponding elements in
$\ccliff{p+q,p+q}$, respectively).

Let $G_0\in \KO_0$ be an operator of graded index $1$ (this is an
operator of type $Inf$ as described in the corresponding $KO$ and
$ko$ models in \cite{HST, markert}; compare also with
\cite[section 3, p. 302]{joachimKO}, where the Fredholm operator
models are used for this construction). Then we take $G_{0, 0}$ to
be the image of $G_0$ via the natural map $\ko_0= \KO _0 \to
\KR_{0, 0}= \kr_{0,0}$. We could start directly  with any $G_{0,0}
\in \KR_{0,0}$ of graded index $1$ in the fixed points, but by
choosing $G_{0,0}$ coming from $\KO_0$, we ensure that the natural
map $\KO_\bullet \to \KR_{\bullet, 0}$ is a map of symmetric
spectra.

\sm

We now define a map
$$
\RR ^{p,q} \rightarrow \kr_{p,q},~~~ v\mapsto L_v\star G_{0,0}^{\star p+q}~~.
$$
and a map
$$
\RR ^{p,q} \rightarrow \KR_{p,q}
$$
by the commutative diagram of $\Sigma_p \times
\Sigma_q$-equivariant maps:
\begin{eqnarray}
\label{pre:unit:equiv:1}
\xymatrix{
\RR ^{p,q} \ar[r]^{\eta_{\kr}}\ar[rd]^{\eta_{\KR}}&\kr_{p,q} \ar[d]^{\covermap_{\KR}}\\
&\KR_{p,q} }
\end{eqnarray}

\sm These extend to continuous \emph{pointed} maps
$$
\eta_\kr:\overline{\RR} ^{p,q}= \sphere{p,q}\rightarrow
\kr_{p,q},  ~~\eta_\KR: \overline{\RR}^{p, q}= \sphere{p,
q}\rightarrow   \KR_{p, q}
$$
by sending base points to base points.
$\overline{\RR}^{p,q}=\sphere{p,q}$ denotes the one point
compactification of the locally compact space $\RR^{p,q}$.

In other words, Diagram (\ref{pre:unit:equiv:1}) extends to a
commutative diagram of $\Sigma_p \times \Sigma_q$-equivariant maps
\begin{eqnarray}
\label{unit:equiv}
\xymatrix{
\sphere{p,q} \ar[r]^{\eta_{\kr}}\ar[rd]^{\eta_{\KR}}&\kr_{p,q} \ar[d]^{\covermap_{\KR}}\\
&\KR_{p,q} }
\end{eqnarray}

\sm

We draw the reader's attention on a small simplification allowed by considering operators of $\hinf$ type instead of for example Fredholm operators as in \cite{joachimKO}: we do not need to extend our space of infinitesimal operators to get the unit.


\subsubsection{Equivariant ring spectrum structure}

\begin{theorem}
The  symmetric sequences of spaces $\kr=(\kr_{n,n})_{n\geq 0}$ and
$\KR=(\KR_{n,n})_{n\geq 0}$ with the above symmetric action,
multiplication and units form commutative $\ztwo$-equivariant
symmetric ring spectra. $\KR$ is a representing spectrum for
$\K$-theory with reality. The map $\covermap_{\KR}$ is an
equivariant ring map.
\end{theorem}

We implicitly assert that the spectrum structure of $\KR$ obtained
from the symmetric structure is exactly the one that corresponds
to Bott periodicity (1-1 periodicity, actually), and that the
product structure lifts that in $K$-groups. This does not exclude
the possibility of the existence of different ring structure
\emph{at the spectrum level} although one expects uniqueness as
highly structured equivariant ring spectrum in view of \cite{BR1}.
In particular, the adjoints of the structure maps for the
symmetric $\ztwo$-spectrum structure of $\KR$
$$
b_{0,1}: \KR_{p, q} \to \KR_{p, q+1}, ~~ b_{1,0}: \KR_{p, q} \to
\KR_{p+1, q},~~ b_{1,1}:= b_{1,0}b_{0,1}: \KR_{p, q}\to \KR_{p+1,
q+1}$$ are \emph{equivariant} weak homotopy equivalences for all
$p,q\geq 0$ (that is also after taking $\ztwo$-fixed points).
Moreover, the diagrams

{\small
\begin{eqnarray}
\label{diagram:comm:bott} \xymatrix{\ar[d]^{b_{0,1}} \kr_{p, q}
\ar[r]^{\covermap}& \ar[d]^{b_{0,1}}\KR_{p,q} &
\ar[d]^{b_{1,0}} \kr_{p, q} \ar[r]^{\covermap}& \ar[d]^{b_{1,0}}\KR_{p,q}\\
\Omega^{0,1}  \kr_{p,
q+1}\ar[r]^{\Omega^{0,1}\covermap}&\Omega^{0,1}  \KR_{p,q+1}&
\Omega^{1,1}  \kr_{p+1,q}\ar[r]^{\Omega^{1,0}\covermap}&\Omega^{1,0}  \KR_{p+1,q}\\
&\ar[d]^{b_{1,1}} \kr_{p, q} \ar[r]^{\covermap}& \ar[d]^{b_{1,1}}\KR_{p,q} \\
&\Omega^{1,1}
\kr_{p+1,q+1}\ar[r]^{\Omega^{1,1}\covermap}&\Omega^{1,1}
\KR_{p+1,q+1}}
\end{eqnarray}}
commute for $p, q\geq0$.

\sm The fact that our spectrum represents usual $\KR$-theory
follows from a sequence of observations. Beginning with the
classical fact that the space $\KR_{0,0}$ is a representing space
for Atiyah's $\KR$ functor on Real spaces (that is, spaces with
$\ztwo$-action), we observe that the spaces $\KR_{p,q}$ are
representing spaces for certain bundles with Clifford structure on
compact spaces.

From \cite[section 5]{atiyahsinger} we see that the multiplication
given by the symmetric spectrum structure gives back the tensor
product of bundles. Finally $1-1$-periodicity corresponds to
tensoring with a generator of $\KR_{1,1}$. This means that the
spectrum structure and the multiplication are the right one in the
homotopy category. We implicitly assume the natural homotopy
equivalence between the spaces of operators of $\hinf$ type and
the (ungraded version of) operators of Fredholm type. The way to
compare our model with the classical models is thoroughly
explained in \cite{HST} and \cite[Section 3.1]{markert}.

\section{Connective $\K$-theory with reality}

\label{connective:dugger}

In \cite{dugger} the author develops a new
proof that the natural spectrum map $KU^{\ztwo}\to KU^{h\ztwo}$ is
a weak equivalence. To this end he introduces a spectrum $\kr$, a
connected version of $\KR$, which is well behaved under taking
fixed points (\cite[section 6, proof of Proposition 6.1]{dugger}). This produces the
following characterization of $\kr = \{\kr_{n,n}\}_{n\geq 0}$.

\begin{proposition}

\label{charac:dugger} There is up to levelwise $\ztwo$-equivariant
weak homotopy equivalence a unique $\ztwo$-equivariant spectrum
$\kr= \{\kr _{n, n}, b_{1,1}\}_{n\geq 0}$ together with a
$\ztwo$-equivariant map of spectra $\{(\covermap) _{n,n}:
\kr_{n,n} \rightarrow \KR_{n,n}\}$ such that:

\begin{enumerate}

\item[(i)]\label{point1} $\kr_{n,n}$ is $(2n-1)$-connected,
$[\kr_{n,n}]^{\ztwo}$ is $(n-1)$-connected,

\item[(ii)]\label{point2}  $\pi_i((\covermap) _{n,n}):
\pi_i(\kr_{n,n} )\to\pi_i(\KR_{n,n} )$ is an isomorphism for
$i\geq 2n$ and $\pi_i([\covermap
]^{\ztwo}):\pi_i([\kr_{n,n}]^{\ztwo}) \to\pi_i([\KR_{n,n}
]^{\ztwo}) $ is an isomorphism for $i\geq n$.

\end{enumerate}

\end{proposition}

Here the superscript $[-]^{\ztwo}$ means fixed points under the $\ztwo$-action.

In fact, any spectrum satisfying the properties listed above has a
$\ztwo$-equivariant map to D. Dugger's $\kr$, which one checks to
be a $\ztwo$-equivariant weak homotopy equivalence. It is useful
to observe here that a $\ztwo$-equivariant spectrum indexed on a
complete universe is determined by its `diagonal spectrum'.

\sm One has to keep in mind that in \cite{dugger} the author
chooses to use a \emph{motivic indexing} that differs from ours.
For instance, the $\Omega ^{2,1}$ of \cite{dugger} corresponds to
our $\Omega ^{1,1}$. We show

\begin{theorem}
The spaces $kr_{p,q}:=\hinnfR{p,q}{q,p}$ form a spectrum which is
levelwise $\ztwo$-equivariantly weakly equivalent to the spectrum
$kr$ of \cite{dugger}.
\end{theorem}

We prove this statement by showing that our construction
$\covermap: \kr \to \KR$ satisfies the properties in the
characterization in Proposition \ref{charac:dugger}. We now give
an overview of the proof.

\subsection{Overview of the method}

\label{overview}

The crucial technical point is the construction of
$\ztwo$-equivariant weak homotopy equivalences (in Section
\ref{quasi:defs}):
\[
\bottonezero: \hinnfR{p,q}{q,p}\longrightarrow
\Omega^{1,0}\hinnfR{p,q+1}{q+1,p}
\]
\[
\bottzeroone: \hinnfR{p,q}{q,p}\longrightarrow
\Omega^{0,1}\hinnfR{p+1,q}{q,p+1}.
\]
This means that the maps induce weak homotopy equivalences on
fixed points as well. The maps will be produced from certain
equivariant quasi-fibrations with equivariantly contractible total
space.

\sm Next we observe that the space $\hinnfR{p, q}{k,l}$  is
connected (Section \ref{connectivity}) for $p,q,k,l\geq 0$ and
$k+l\geq 1$, as well as its fixed point space $[\hinnfR{p,
q}{k,l}]^{\ztwo}$ for $k\geq 1$. This is done by finding explicit
paths to the base point (through fixed points if the starting
point was a fixed point). Hence by the above equivalences, the
space $\kr _ {p, q}$ is at least $(p+q-1)$-connected. On fixed
points, we find that $[\kr_{n, n}]^{\ztwo}$ is $(n-1)$-connected.
This settles point (\ref{point1})$(i)$ in the characterization of
$\kr$.

\sm The  point $(ii)$ in Proposition \ref{point2} is slightly more
involved, and relies on further properties of the quasifibrations.
We let

$$
\bottoneone := \bottonezero\bottzeroone:
\kr_{p,q}=\hinnfR{p,q}{q,p}\longrightarrow
\Omega^{1,1}\hinnfR{p+1,q+1}{q+1,p+1}= \Omega^{1,1}\kr_{p+1, q+1}
$$

In Section \ref{commutation}, we show that the following diagram  of $\ztwo$-equivariant maps  commutes

\begin{eqnarray}
\label{comm} \xymatrix{
\ar[d]^{\bottoneone}\kr_{n, n} \ar[r]_{=}& \ar[d]^{b_{1,1}} \kr_{n, n} \ar[rr]^{\covermap _{n,n}}&& \ar[d]^{b_{1,1}}\KR_{n,n} \\
\Omega^{1,1} \kr_{n+1,n+1} \ar[r]_{=}&\Omega^{1,1}
\kr_{n+1,n+1}\ar[rr]_{\Omega^{1,1}\covermap _{n+1,
n+1}}&&\Omega^{1,1}  \KR_{n+1,n+1}~. }
\end{eqnarray}
Under forgetting the $\ztwo$-action, we get inductively, beginning
with $\kr_{0,0}= \KR_{0,0}$, that $\pi_k(\covermap_{n,n}):
\pi_k(\kr_{n, n})\to \pi_k(\KR_{n,n})$ is an isomorphism for
$k\geq 2n$. Indeed, the commutation of Diagram (\ref{comm})
together with the fact that $\bottoneone$ and $b_{1,1}$ are
$\ztwo$-equivariant  weak equivalences, yields the following. If
$\pi_k(\covermap_{n,n}):\pi_k  (\kr_{n,n}) \to  \pi_k(\KR_{n,n})$
for $k\geq 2n$,  as $\pi_k (\Omega^{1,1}\covermap_{n+1,n+1 })$ and
$\pi_k(\covermap_{n,n})$ are  isomorphisms in the same dimensions
$k$, and  as  $\Omega^2= \Omega ^{1,1}$ by forgetting the action,
we get that
$$
\pi_k(\covermap):\pi_k(\kr_{n+1, n+1}) \to\pi_k( \KR_{n+1, n+1})
$$
is an isomorphism for $k\geq 2n+2$.

\sm Taking fixed points in Diagram (\ref{comm}) yields by a  similar but slightly more complicated argument that
$$\pi_k ([\covermap]^{\ztwo}): \pi_k( [\kr_{n,n}]^{\ztwo} )\to\pi_k ([\KR_{n,n,}]^{\ztwo}) $$
is an equivalence in degrees $k\geq n$, which is the last
statement to be proved in the characterization. The details of
this argument are deferred to Section \ref{commutation}.

\section{Proofs}

\label{proofs}

\subsection{CW-structures}

\label{CW}

It is an important if somewhat technical point to check that our spaces have the homotopy type
of CW-complexes with a nice filtration by subcomplexes.

\sm This property is needed later on to allow us to use the
Dold-Thom criterion for quasi-fibrations. It was shown in
\cite{markert} in the case of the spaces for connective real
$\K$-theory $\ko$ building on a similar proof in \cite[section
7]{HST}. The approach works under minor modifications. We give a
short outline here, in order to enable the reader to trace the
differences in the $\ztwo$-equivariant context. The filtration is
constructed as follows: let
$$
\ldots \subset \hilb{}{}(i) \subset \hilb{}{}(i+1)\subset \ldots
\subset \hilb{}{}
$$
for $i\geq 0$ be a filtration of $\hilb{}{}$ by finite dimensional
subspaces of increasing dimension. This induces a filtration
$$
\ldots \subset \hilbR{p, q}{k,l}(i) \subset \hilbR{p,
q}{k,l}(i+1)\subset \ldots \subset  \hilbR{p, q}{k,l}
$$
by graded free $\ccliff{q,p}-\ccliff{k,l}$-submodules, for any
$p,q,k,l\geq 0$. In the following, this indexing does not change
and so we omit the indices of the Clifford-structures. We obtain a
filtration
$$
\ldots \subset \hinnfR{}{}(i)\subset \hinnfR{}{}(i+1)\subset
\ldots \subset \hinnfR{}{}
$$
where $\hinnfR{}{}(i):= \hinnfR{}{}(\hilbR{}{}(i))$. Recall that
the definition of the $\hinnfR{}{}$ spaces also makes sense in the
context of a finite dimensional Hilbert universe. Observe first
that the natural map
$$
\colim_i ~([\hinnfR{}{}(i)]^{\ztwo} )\to [\colim_i
(\hinnfR{}{}(i))]^{\ztwo}
$$
is a homeomorphism, because the colimit is taken over
\emph{monomorphisms}.

\begin{proposition}\label{proposition:CW}

The spaces $\hinnfR{}{}(i)$ as well as their fixed point spaces
$[\hinnfR{}{}(i)]^{\ztwo}$  have the homotopy type of CW
complexes. The natural inclusions
$$
\colim_i~ (\hinnfR{}{}(i)) \to \hinnfR{}{}, \quad\quad\quad
\colim_i( [\hinnfR{}{}(i)]^{\ztwo}) \to[\hinnfR{}{}]^{\ztwo}
$$
are homotopy equivalences.
\end{proposition}

In particular, $\colim_i~ \hinnfR{}{}(i) \simeq \hinnfR{}{}$ and $\colim_i ~([\hinnfR{}{}(i)]^{\ztwo}) \simeq [\hinnfR{}{}]^{\ztwo}$ have both the homotopy type of CW complexes.

\sm

{\it Proof of Proposition \ref{proposition:CW}.~} The fact that
the inclusions of the colimits are homotopy equivalences is an
{\it ad hoc} adaptation of the arguments in \cite[Proposition 29
and Remark 30]{HST}.

\sm

The fact that the spaces are CW-complexes up to homotopy is more
involved. The spaces $\hinnfR{}{}(i)$ are homeomorphic to
realizations of nerves of certain topological categories
$\mathfrak{C}(i)$ (of course these categories should be indexed by
the appropriate Clifford structure, i.e.
$\mathfrak{C}_{p,q}^{k,l}(i)$, we suppress this again).

\[
\centerline{
\xymatrix{ \ldots\ar[r] & \hinnfR{}{}(i)\ar[r]  & \hinnfR{}{}(i+1)
\ar[r] & \ldots\ar[r] & \colim_i \hinnfR{}{}(i)\ar[r]^{~~\sim} & \hinnfR{}{}\\
\ldots\ar[r] & |\nerve (i)|\ar[r]\ar[u]^{\homeo} &  |\nerve
(i+1)|\ar[r]\ar[u]^{\homeo} & \ldots\ar[r] &  \colim_i  |\nerve
(i)|\ar[u]^{\homeo}
& \\
}}
\]

The nerve of these categories are simplicial topological spaces
$m\mapsto \nerve (i)[m]$ which are \emph{good} in Segal's sense,
and are levelwise CW-complexes. Goodness in this sense implies
that the realization of $\nerve (-)$ is homotopy equivalent to the
fat realization (see \cite{segal}), which itself preserves
$CW$-complex structure. Since we can show this structure
levelwise, we obtain $CW$-complexes in the lower line. The maps
in the bottom colimit system come from the realization  of
levelwise inclusions of $CW$-complexes.

\sm

We define the small topological categories  $\mathfrak{C}(i)$
analogous to the \emph{internal space category} of \cite{HST}. The
objects of $\mathfrak{C}(i)$ are Clifford-linear (w.r.t. the full
$\ccliff{p,q}-\ccliff{k,l}$ structure) subspaces of
$\hilbR{p,q}{k,l}(i)$, topologized as space of finite-rank,
Clifford-linear orthogonal projection operators. The set of
morphisms between $W_0$ and $W_1$ is empty except if $W_0\subseteq
W_1$, in which case its elements are odd, orthogonal involutions
on $W_1 - W_0 := W_0^{\perp_{W_1}}$ which preserve the
\emph{right} Clifford action ($\ccliff{p,q}$-linear), that is
{\footnotesize
\[
\mathfrak{C}(i)(W_0, W_1) = \set{ R \in \OO^{odd}(W_1 - W_0) | R^2
= Id, Re_i=e_iR, i=1, \ldots, p,Rf_j=f_jR, j=1, \ldots, q }.
\]}

The morphism sets are again topologized as subspaces of bounded
operators on $\hilbR{}{}(i)$. The action of $\theta$ on the
$\hinnfR{}{}(i)$ induces an action on object and morphism spaces
of $\mathfrak{C}(i)$.

For the fixed-point case the categories consist of
$\ztwo$-invariant subspaces together with involutions which
commute with $\theta$.

\sm The nerve of $\mathfrak{C}(i)$ is the usual simplicial
topological space made of chains of composable morphisms, with the
product topology.





In particular, a point in $\nerve (i)[m]$ is a length-$m$ flag $(W_0 \subset \ldots\subset  W_m)$ of Clifford-modules together with Clifford-right-linear involutions on the steps $W_i-W_{i-1}$.

\sm The realization of this simplicial space is the  quotient
$$
\coprod_{m\geq 0} \nerve (i)[m]\times \Delta^m / \sim
$$

by the standard relation $\sim$ \cite{HST}. We think of $\Delta^m$
as $m$-tuple $(t_1 \leq \ldots \leq t_m)$ of points in $\bar{\RR}=
[0, +\infty]$.

Now it is clear how to obtain an operator in $\hinnfR{}{}(i)$ from
a point in the realization: we send $(x; t) \in \nerve
(i)[m]\times \Delta^m$ to the operator $G(x; t)$ with spectral
configuration given by:

\begin{itemize}

\item the eigenvalues $0$, $\pm t_i$ for $i=1, \ldots, m$ (notice that these can be listed with multiplicities, $0 =\lambda_0\leq  \lambda_1 \leq \lambda _2 \leq \ldots \leq \infty$; an eigenspace at infinity for $t_i=\infty$ is formally included in the configuration with eigenspace denoting the complement of the domain)

\item  the eigenspace at $0$ is $W_0$, the eigenspaces for $\pm t_i$ are the $(\pm 1)$-eigenspaces of the involution $R_i$ on the module $W_i-W_{i-1}$ (to the formal eigenspace at infinity we add the complement of $W_m$. Notice that the single eigenspaces are linear w.r.t. the right Clifford action, however the eigenspaces $W_i-W_{i-1}$ of the square $G(x; t)^2$ are linear w.r.t. both actions).




\end{itemize}

This produces a continuous map on the quotient $|\nerve (i)|$
which is a homeomorphism (compare \cite[Chap. 7, lemma 98]{HST}).
Again it is obvious that this construction, when applied to the
smaller categories in the fixed-point case, produces operators in
$[\hinnfR{}{}(i)]^{\ztwo}$ since all eigenspaces will be invariant
under the action.

Now one shows:

\begin{itemize}\label{structure:nerve}

\item[-] for each $m\geq 0$, the space $\nerve (i)[m]$ is a CW-complex,

\item[-] the simplicial space $m\mapsto \nerve (i)[m]$ is \emph{good} in Segal's sense.

\end{itemize}











%






The idea for the first part is that each space $\nerve (i)[m]$ is
the total space of a fiber bundle over a flag space, where both
base and fiber have the homotopy type of $CW$-complexes. Namely,
there is an obvious map $p_m(i): \nerve
(i)[m]\longrightarrow\flag_m(i)$  taking an element in $\nerve
(i)[m]$ to its flag $(W_0 \subset \ldots\subset  W_m)$.

The flag space $\flag_m(i) = \coprod_{d=(d_0 \leq d_1\leq
\ldots\leq d_n)\in\NN^{m+1}} \flag_m^d(i)$ is the space of
Clifford-linear flags of length $m$ in $\hinnfR{}{}(i)$. It
consists of components determined by the Clifford dimension of the
successive subspaces of a flag.

Now the map $p_m(i)$ restricts to a fiber bundle over each
component $\flag_m^d(i)$, where the fiber is the space of odd,
orthogonal, $\cliff{p, q}$-linear involutions on the flags (i.e.
involutions on $W_m-W_0$ which keep the subspaces $W_i-W_{i-1}$ of
fixed dimensions invariant).

The flag spaces $\flag_m^d(i)$ are $CW$ since they are themselves total
spaces of bundles over the $d_m$-dimensional Clifford-linear
Grassmannian of $\hinnfR{}{}(i)$, which is $CW$, with $CW$
fibers. The subspace of involutions in
$\mathcal{O}^{odd}_{\ccliff{p, q}}(W_m-W_0)$ preserving a specific
flag in $\flag_m^d(i)$ also has a $CW$-complex decomposition.

\sm Furthermore, the simplicial space $\nerve (i)$ is \emph{good},
which means that the inclusions of the degenerated subspaces

\[
\nerve (i)[m,k]:= s_k(\nerve (i)[m-1])\hookrightarrow \nerve
(i)[m]
\]

are strict neighbourhood retracts in $\nerve (i)[m]$. This can be
easily seen by observing that the spaces $\nerve (i)[m,k]$ map
\emph{surjectively onto connected components} of $\nerve (i)[m]$,
namely those which lie over connected components $\flag_m^d(i)$ of
$\flag _m(i)$ with dimension vectors $d$ of the form
$$
d= (d_0, \ldots, d_k, d_k, \ldots, d_{m-1}).
$$

\sm Using \cite{PS}, we observe that all spaces which  we need to be CW-complex are defined by algebraic equations and are actually $\ztwo$-semialgebraic sets, hence possess a $\ztwo$-CW complex structure. In particular, under taking $\ztwo$-fixed points, the whole argument works again, and this yields the result stated in Proposition \ref{proposition:CW}.





%




\subsection{Constructing the structure maps $\bottonezero$ and $\bottzeroone$.}



Here we define two maps
\[
\phi: \HEE{p,q}{k,l}\longrightarrow\hinnfR{p-1,q}{k+1,l}
\quad\quad\quad\quad \psi:
\HFF{p,q}{k,l}\longrightarrow\hinnfR{p,q-1}{k,l+1}
\]
which we prove to be quasi-fibrations with the same fibers
$\hinnfR{p,q}{k,l}$ and with contractible total spaces. This
implies weak homotopy equivalences (see explicit Formulas
(\ref{old1}) and (\ref{old2}))
\[
\bottonezero: \hinnfR{p,q}{k,l}\longrightarrow
\Omega^{1,0}\hinnfR{p-1,q}{k+1,l} \quad\quad\quad\quad
\bottzeroone: \hinnfR{p,q}{k,l}\longrightarrow
\Omega^{0,1}\hinnfR{p,q-1}{k,l+1}.
\]

In a non-equivariant context, one can compare a quasifibration
with contractible total space with the
path-loop fibration of the base space; this yields two long exact
sequences in homotopy and one obtains weak homotopy equivalences
as above by the five-lemma, since both total spaces are
contractible. In the $\ztwo$-equivariant case, we have to take
care that the actions are consistent.  This leads to the two
different versions of the loop space.

\sm We can reformulate the maps $\bottonezero$ and $\bottzeroone$
to adjust the indexing.

This is done by using $(1,1)$-periodicity first, and then applying the maps $\phi$ or $\psi$, respectively.

This produces the weak equivalences
$$\bottonezero: \hinnfR{p,q}{k,l}\longrightarrow\Omega^{1,0}\hinnfR{p,q+1}{k+1,l}\quad \quad\quad\quad
\bottzeroone:
\hinnfR{p,q}{k,l}\longrightarrow\Omega^{0,1}\hinnfR{p+1,q}{k,l+1}.$$

For $\kr_{p,q}$, this yields the deloopings we need.

The index shift using $(1,1)$-periodicity can be written out,
analogously to \cite[section 5]{joachimKO}, as the equivariant
homotopy equivalence
\begin{eqnarray}\label{formula:shift}
\hinnfR{p,q}{k,l}\simeq \hinnfR{p+1,q+1}{k,l}, ~ G\mapsto
(Id_{\ccliff{1,1}}\gradedtensor G_{0,0})\star G
\end{eqnarray}
where $G_{0,0}$ is the generator of $\KR_{0,0}$ as before. Note
that $\widetilde{G_{0,0}}:= Id_{\ccliff{1,1}}\gradedtensor
G_{0,0}$ is defined on $\hilbR{1,1}{}$, therefore the image above
is an operator on $\hilbR{p,q}{k,l}\gradedtensor
\hilbR{1,1}{}\cong \hilbR{p+1,q+1}{k,l}$ as wished.
The fact that the above map is a homotopy equivalence is fully
explained in \cite{joachimKO} and relies in turn heavily on the
seminal paper \cite{atiyahsinger}.

\sm




To prove that the maps a quasifibration, we use the following version of a well-known theorem of Dold and Thom \cite{DT}:




\begin{thm}\label{thm:doldthom}\label{Dold-Thom:criterion} Let $q:X\rightarrow Y$ be surjective. Then $q$ is a quasifibration, if there is an increasing filtration $\set{F_iY}$ of $Y$, where $\colim_{i}F_iY\sim Y$ is a homotopy equivalence, such that the following holds:

\begin{itemize}

\item[(i)]  For every open subset $U$ of $B_i:=F_i-F_{i-1}$, the restriction $\pminusone{U} \longrightarrow U$ is a fibration.

\item[(ii)] For every $i$, there exist neighbourhoods $N_i \subset F_{i+1}$ of $F_i$ and a contracting homotopy $h: N_i\times I \rightarrow N_i$ with $h_0=id_{N_i}$ and $h_1(N_i)\subset F_i$.

\item[(iii)] This deformation $h$ is covered by a homotopy $H: q^{-1}(N_i)\times I \rightarrow q^{-1}(N_i)$, $H_0 = id_{q^{-1}(N_i)}$,  such that for each point $x \in N_i$

\[H_1: q^{-1}(x) \rightarrow q^{-1}(h_1(x))\] is a weak homotopy equivalence.

\end{itemize}

\end{thm}

\sm

In the following, given an increasing filtration $\{X_i\}_{i\geq
0}$ of a space $X$, the spaces $X_{i+1} -X_{i}$ are called the
filtration strata  and the space $X_i$ will be called the $i$th
filtration step.

\sm

It is not difficult to show that the maps $\phi$ and $\psi$ are
surjective and have the correct fiber over the base point. As we
will see below, there is also an obvious way to filter the base
space, namely by filtering the operators by the Clifford dimension
of their domains. If we think of the operators as given by their
spectral configuration, this means that we filter such
configurations by the  Clifford dimension of the algebraic sum of
eigenspaces. Neighborhoods of one such filtration level within
the next larger one are given by operators whose spectrum has at
least one pair of eigenvalues within a fixed contractible
neighborhood of the base point $\infty$ of $\Rinfty$. The
homotopies $h$ contracting these neighbourhoods are induced by
contracting the neighborhood of $\infty$: we apply functional
calculus with such a homotopy to the operator space, this is
continuous and produces a homotopy as wished. It is also easy to
see that this can be covered by a corresponding homotopy $H$ in
the total spaces. Thus the constructive parts of this criterion in
Theorem \ref{Dold-Thom:criterion} are not difficult (including the
point ii) of this criterion); it remains to show that the maps are
fibration on the filtration strata (i.e. on the subspaces of
operators with fixed Clifford dimension of domain), and that the
homotopies contracting neighborhoods of filtration levels induce
weak homotopy equivalences in the fibers. Indeed the fibers over
different filtration strata are not homeomorphic and we need an
elaborate technical argument to show this last point. We show the isomorphism on
homotopy groups by studying maps from spheres into the fiber
spaces in question.

\sm

Details  of this proof are explained in \cite{markert}, where it
is done for the ordinary real connective theory $\ko$, and we will
omit them.

In the  next section we will provide the definitions of the total
spaces and maps and the essential parts of the proofs with
particular respect to the adjustments in the $\ztwo$-equivariant
context.

\subsubsection{Definition of the quasi-fibration maps.}

\label{quasi:defs}

We define an operator class $\hinnfR{p,q}{ev, -\infty}$ of {\it
even}, $\ccliff{p,q}$-linear operators of type $\hinf$, with
eigenvalues of finite multiplicity in $[-\infty, +\infty)$. This
means that we can write the operator as a formal sum $E = \sum
\lambda\pi_{E_\lambda}$ where $\lambda \in [-\infty, +\infty)$.
The complex equivalent of this is a class
$\hinnfR{p,q}{ev,-i\infty}$ of {\it even, skew-adjoint}
Clifford-linear operators with finite-rank resolvent and spectrum
in $[- i\infty, + i\infty)$.

\sm

As a set, the total space $\HEE{p,q}{k,l}$ consists of triples
$(\D; G,E)$ where $\D\subset\hilbR{p,q}{k,l}$ is a graded,
$\ccliff{p+k,q+l}$-linear closed subspace (the domain); $G$ is an
operator in $\hinnfR{p,q}{k,l}$ and $E$ an operator in
$\hinnfR{p+k,q+l}{ev, -\infty}$ such that $\dom{G}=\dom{E}=\D$ and
$[G,E]=0$ on $\D$. The set $\HFF{p,q}{k,l}$ is defined analogously
as triples $(\D; G, F)$ of commuting operators
$G\in\hinnfR{p,q}{k,l}$, $F\in\hinnfR{p+k,q+l}{ev,- i\infty}$,
defined on a common domain $\D$.

\sm

Just as we can think of an operator $G\in\hinnfR{p,q}{k,l}$ as a
configuration of points in $\Rinfty$ (eigenvalues) with
coefficients (projections onto eigenspaces), we can think of a
triple $(\D; G,E)\in\HEE{p,q}{k,l}$ as a configuration of points
on the cone $[-\infty, +\infty]\wedge\Rinfty$ with coefficients as
before (the base point of $[-\infty, +\infty]$ is at $+\infty$).
The points and coefficients now correspond to the eigenvalues and
eigenspaces in a {\it simultaneous} spectral decomposition of the
two commuting operators $G$ and $E$. In the second case the same
holds true if we replace the cone by a complex cone $[-i\infty,
i\infty]\wedge\Rinfty$. This small difference will be crucial for
the commutativity with the non-trivial $\z2$-action on the loop
space $\Omega^{0,1}$ (compare \cite[Section 2.2]{schro}).

\sm

The topology of the total spaces is defined in terms of these
configuration spaces (compare remarks following the definition of
the spaces $\hinnfR{p,q}{}$). To put a precise expression in the
notation\footnote{The superscript $\ztwo$ from \cite{markert} in the notation for
the configuration spaces refers to a $\ztwo$-grading.
To makes things clear, we use the notation $[-]^{\ztwo}$ to
express fixed point spaces in this work. } of \cite{markert}, we
want to have
$$
\HEE{p,q}{k,l}:=\Conflinfin{\ccliff{p+k,q+l}}{\z2}([-\infty,
+\infty]\wedge\Rinfty;\hilbR{p,q}{k,l})\quad\quad\quad\quad
\HFF{p,q}{k,l}:=\Conflinfin{\ccliff{p,q}}{\z2}([-i\infty,
i\infty]\wedge\Rinfty;\hilbR{p,q}{k,l}).
$$

The sets of triples and configurations as above are homeomorphic.
For our calculations we will use the operator description. Both
total spaces are $\ztwo$-spaces by the involution $\theta: (-,
-)\mapsto (\theta - \theta,\theta - \theta)$.

\sm Now the maps are defined as

\begin{eqnarray}
\label{formula1} \phi:
\HEE{p,q}{k,l}\longrightarrow\hinnfR{p-1,q}{k+1,l}: (\D; G,
E)\mapsto (G+ E\epsilon e_p)
\end{eqnarray}
\begin{eqnarray}
\label{formula2} \psi:
\HFF{p,q}{k,l}\longrightarrow\hinnfR{p,q-1}{k,l+1}: (\D; G,
F)\mapsto (G + F \epsilon f_q)
\end{eqnarray}
where the image is defined on $\D\cap\Esp{E}{-\infty} ^{\perp}$
(or $\D \cap \Esp{F}{-i\infty}^{\perp}$ respectively). Here $e_p$ is
the $p$th generator of $\ccliff{p,q}$ and we think of $\epsilon
e_p$ as the $(k+1)$st generator of $\ccliff{l, k+1}$ acting from
the left, while $f_q$ is the $q$th generator of $\ccliff{p,q}$ and
we think of $\epsilon f_q$ as the $(l+1)$st generator of
$\ccliff{l+1, k}$ acting from the left.

\sm

Note that an image operator $\phi(G,E)$ satisfies $(\phi(G,E))^2 =
G^2 + E^2$, while for an element $\psi(G,F)$ in the image of
$\psi$ we have $(\psi(G,F))^2 = G^2 - F^2$. This means that the
operator in the image of either of the maps has positive real
eigenvalues of the shape $\pm\sqrt{\lambda^2 + \mu^2}$ for
eigenvalues $\lambda$ of $G$ and $\mu$ of $E$ ($\lambda$ of $G$
and $i\mu$ of $F$, respectively, in the case of $\psi$) and the
corresponding eigenspaces are the intersections of eigenspaces of
$G$ and $E$ (resp. $F$) at these eigenvalues.

\sm

One can determine the fibers over a point $P$ in the base space of either fibration

$$\phi^{-1}(P) \homeo \hinnfR{p,q}{k,l}(\domperp{P}) \quad \quad \psi^{-1}(P) \homeo \hinnfR{p,q}{k,l}(\domperp{P})$$
by calculating the graded commutator of $P$ with $\epsilon e_p$
(or $\epsilon f_q$ respectively), which recaptures the information
of the pre-image operators on the part $\dom{P}$ of their domain.

\sm Furthermore, as configuration spaces on cones, the spaces
$\HEE{p,q}{k,l}$ and $\HFF{p,q}{k,l}$ are contractible by applying
functional calculus with a suitable contraction of the cones to
the operators. This can be done $\ztwo$-equivariantly, since
functional calculus affects only eigenvalues, not eigenspaces, and
therefore commutes with $\theta$.

\sm We now check the part i) and iii) of the Dold-Thom criterion (Theorem \ref{Dold-Thom:criterion}) for our maps.

\subsubsection{Proof of part $(i)$: fibration on filtration strata}



We use filtration levels
\[F_i =  \left\{\begin{array}{l|l}P\in\hinnfR{p-1,q}{k+1,l} & \mathrm{dim}_{\cliff{p+k,q+l}}\dom{P}\leq i \\\end{array}\right\}\]
and filtration strata
\[B_i =  \left\{\begin{array}{l|l}P\in\hinnfR{p-1,q}{k+1,l} & \mathrm{dim}_{\cliff{p+k,q+l}}\dom{P}=i \\\end{array}\right\}\]
in the first case and analogously in the second (with different
indexing of the algebras). To satisfy \ref{thm:doldthom}$(i)$ we
have to show that the maps are fibrations on the filtration
strata.



\sm

Within a filtration stratum $B_i$ the domains of the operators
form a subspace of a finite Grassmannian
$\fingrass{\hilbR{p,q}{k,l}}{i}$ of $\ccliff{p+k,q+l}$-linear
subspaces of $\hilbR{p,q}{k,l}$. The pullback of the complement of
the tautological bundle over this Grassmannian determines the
subspaces $\domperp{-}$ which distinguish the fibers.

The local trivializations of this bundle induce local trivializations of both $\phi$ and $\psi$ on the filtration stratum $B_i$ by conjugating the operators.

\subsubsection{Proof of part $(iii)$: weak equivalence of fibers under change of filtration level}

As discussed above we now need to prove that a deformation $h$
which contracts a neighbourhood $N_i$ of $F_i$ into $F_{i-1}$ is
covered by a homotopy $H$ such that for each operator $P \in N_i$,
the change in the fibers of $\phi$ ($\psi$ respectively)
\[H_1:\phi^{-1}{P}\rightarrow \phi^{-1}{h_1(P)}\] is a weak homotopy equivalence.

\sm We do this for $\phi$, the proof works analogously for the map
$\psi$. The crucial difference in the fibers $\phi^{-1}{P}$ and
$\phi^{-1}{h_1(P)}$ is that the domains of $P$ and $h_1(P)$ differ
by some finite-dimensional $\ccliff{p+k,q+l}$-linear subspace $V$
(the homotopy $h$ removes those summands in the spectral
decomposition of $P$ whose eigenvalues are within a given
neighbourhood of $\infty$).

Now we can show in general that the following holds:

\begin{prop}\label{prop:fiberjump}
An inclusion $i: R\hookrightarrow R\oplus V$ of
infinite-dimensional graded Real
$\cliff{l,k}-\cliff{p,q}$-bimodules with finite codimension
induces a weak homotopy equivalence
\[i_*: \hinnfR{p,q}{k,l}(R) \stackrel{\sim}{\longrightarrow} \hinnfR{p,q}{k,l}(R\oplus V).\]
\end{prop}

This implies in particular that there is a weak homotopy equivalence of the fibers

\[\phi^{-1}{h_1(P)}\stackrel{\sim}{\longrightarrow} \phi^{-1}{P}\]

induced by the inclusion $\domperp{h_1(P)}\hookrightarrow \domperp{P}$.


\sm The proof of proposition \ref{prop:fiberjump} is rather
technical and uses the fact that the involved spaces have the
homotopy type of $CW$-complexes as shown in proposition
\ref{proposition:CW}.

\sm

As infinite-dimensional Clifford-modules the two label spaces $R$
and $R\oplus V$ are isomorphic. Choosing such an isomorphism
$\tau$ gives a homeomorphism and a diagram:

\[
\centerline{\xymatrix{ \hinnfR{p,q}{k,l}(R)\ar@<1ex>[rr]^{i_*} &
&\hinnfR{p,q}{k,l}(R\oplus V)\ar@<1ex>[ll]^{\tau_*}.\\}}\]

We claim that these maps induce inverse to each other isomorphisms
on homotopy groups. For any map $c$ that represents an element
 in $\pi_\bullet(\hinnfR{p,q}{k,l}(R))$ we check that it is
homotopic to $\tau_* \circ i_*\circ c$. To see this we first
transform the problem: $c$ factors up to homotopy as shown

\[
\centerline{ \xymatrix{
\hinnfR{p,q}{k,l}(\hilbR{p,q}{k,l}(i)) \ar[rr]^{j_*} & & \hinnfR{p,q}{k,l}(R) \ar@<1ex>[rr]^{i_*} & & \hinnfR{p,q}{k,l}(R\oplus V) \ar@<1ex>[ll]^{\tau_*} \\
 & & \sphere{\bullet} \ar[u]_{c} \ar[ull]^{\tilde{c}} \\
}}
\]
for some $j$, by Proposition \ref{proposition:CW}. We observe that
this is the point where we need our spaces to have a nice
$CW$-filtration.

\sm

This means that we can instead show that the factorization of $c$
is homotopic to the factorization of $\tau_* \circ i_*\circ c$.
This can be solved on the level of Real Hilbert spaces

\[
\centerline{\xymatrix{ \hilbR{p,q}{k,l}(i) \ar[rr]^{j} & & R
\ar@<1ex>[rr]^{i} & & R\oplus V \ar@<1ex>[ll]^{\tau}\\}}
\]

where the maps $j$ and $\tau \circ i \circ j$ are
$\ztwo$-equivariantly homotopic as $\ztwo$-equivariant embeddings
of $\ccliff{p+k, q+l}$-modules of infinite codimension (the space
of such embeddings is equivariantly contractible). This induces
homotopies on the maps of operators (compare general properties of
the operator spaces, \cite{markert}).

\proofbox


Thus we can complete the proof that both maps are indeed (non-equivariantly) quasi-fibrations in the same fashion as in \cite{markert}.





%

\subsubsection{$\ztwo$-equivariance}

\label{quasi:commute}

The quasifibrations given by the Formulas (\ref{formula1}) and (\ref{formula2}) induce weak homotopy equivalences

\begin{eqnarray}
\label{bott:shift:formula:1}\label{old1} \bottonezero:
\hinnfR{p,q}{k,l}\longrightarrow
\Omega^{1,0}\hinnfR{p-1,q}{k+1,l}: G \mapsto \left( t\mapsto
(\bottonezero G)_t:= G + E_t \epsilon e_p, \,\,\, t\in [-\infty,
+\infty] \right),
\end{eqnarray}

\begin{eqnarray}
\label{bott:shift:formula:2}\label{old2} \bottzeroone:
\hinnfR{p,q}{k,l}\longrightarrow
\Omega^{0,1}\hinnfR{p,q-1}{k,l+1}: G \mapsto \left(t\mapsto
(\bottzeroone G)_t:= G + F_t\epsilon f_q, \,\,\, t\in [-\infty,
+\infty] \right),
\end{eqnarray}

where $E_t := t Id_{\dom{G}}$ and  $F_t := t i Id_{\dom{G}}$.



\sm

We check the commutation with the $\ztwo$-actions. The map
$\bottonezero$ is $\z2$-equivariant: we have

\begin{eqnarray}
(\Omega^{1,0}\theta(\bottonezero G))_t & = & \theta  (G + E_t \epsilon e_p) \theta  \\
 & = & \theta G \theta + \theta E_t \theta \epsilon e_p \nonumber\\
 & = & \theta G \theta + E_t\epsilon e_p\nonumber\\
 &=& (\bottonezero(\theta G))_t \nonumber .
\end{eqnarray}
since $E_t$ commutes with $\theta$
\begin{eqnarray}
\theta  E_t \theta  (x\tensor v) & = & \theta  E_t (\conjugate{x}\tensor \theta(v))\\
& =& \theta (t\conjugate{x}\tensor \theta(v))\nonumber \\
& = & t x\tensor v\nonumber \\
&=& E_t (x\tensor v)\nonumber.
\end{eqnarray}


In the case of $\bottzeroone$ we want to see that the diagram commutes

\[\centerline{
\xymatrix{
\hinnfR{p,q}{k,l}\ar[rr]^{\bottzeroone}\ar[d]^{\theta} & & \Omega^{0,1}\hinnfR{p,q-1}{k,l+1}\ar[d]^{\Omega^{0,1}\theta} \\
\hinnfR{p,q}{k,l}\ar[rr]^{\bottzeroone} & &
\Omega^{0,1}\hinnfR{p,q-1}{k,l+1} } }
\]

where now $\Omega^{0,1}\theta(\gamma)(t) = \theta(\gamma(-t))$ for
$\gamma$ a loop in $\Omega^{0,1}\hinnf{p,q-1}{k,l+1}$, $t\in
[-\infty, +\infty]$. We have

\begin{eqnarray}
(\Omega^{0,1}\theta(\bottzeroone G))_t &=& \theta((\bottzeroone G)_{-t})\\
 & = & \theta (G + F_{-t}\epsilon f_q)\theta  \nonumber \\
 & = & \theta G \theta - \theta F_{-t}\theta\epsilon f_q    \nonumber\\
 & = & \theta G \theta + F_t \epsilon f_q  \nonumber \\
 & =& (\bottzeroone(\theta G))_t  \nonumber.
\end{eqnarray}

The first equation follows since $\theta$ commutes with $\epsilon$
but anti-commutes with the multiplication by $f_q$ (where $f_q$ is
in the $(-1)$-eigenspace of $\theta$). The second identity follows
from
\begin{eqnarray}
 (\theta F_{-t} \theta)(x\tensor v) & = & \theta F_{-t} (\conjugate{x}\tensor \theta(v))  \\
 & = & \theta (-t i \conjugate{x}\tensor \theta(v))  \nonumber\\
& =& \conjugate{-ti\conjugate{x}}\tensor v  \nonumber \\
 & = & tix\tensor v  \nonumber\\
 &=& F_t (x\tensor v)  \nonumber.
\end{eqnarray}
Thus $\psi$ indeed induces a map into the $\Omega^{0,1}$-loop
space.

Finally, combining the maps from Formula
(\ref{bott:shift:formula:1}) and (\ref{bott:shift:formula:2}) with
the $(1,1)$-periodicity given in formula \ref{formula:shift}, we
obtain the (weak) homotopy equivalences which constitute the
structure maps

{\small
\begin{eqnarray}
\label{bott:shift:formula:3} \bottonezero:
\hinnfR{p,q}{k,l}\longrightarrow
\Omega^{1,0}\hinnfR{p,q+1}{k+1,l}: G \mapsto \left( t\mapsto
(\bottonezero G)_t:= \widetilde{G_{0,0}}\star G + \tilde{E_t}
\epsilon e_{p+1}, \,\,\, t\in [-\infty, +\infty] \right),
\end{eqnarray}
\begin{eqnarray}
\label{bott:shift:formula:4} \bottzeroone:
\hinnfR{p,q}{k,l}\longrightarrow
\Omega^{0,1}\hinnfR{p+1,q}{k,l+1}: G \mapsto \left(t\mapsto
(\bottzeroone G)_t:= \widetilde{G_{0,0}}\star G +
\tilde{F_t}\epsilon f_{q+1}, \,\,\, t\in [-\infty, +\infty]
\right),
\end{eqnarray}
where $\tilde{E_t} := t Id_{\dom{ \widetilde{G_{0,0}}\star G}}$
and  $\tilde{F_t} := t i Id_{\dom{ \widetilde{G_{0,0}}\star G}}$.
}

\subsubsection{Quasi-fibration on fixed points}

\label{quasi:fixed}

Finally we want to see that by restricting to fixed-points of the
$\ztwo$-action, the map $\phi$ is also a quasifibration. Recall
that the action on $\hinnfR{}{}$ was given by conjugating with
$\theta$. Fixed points therefore are operators $G$ with $G=\theta
G \theta$: thus they have the same eigenvalues, but
their eigenspaces are closed under the involution $\theta$. Taking
fixed points produces a commutative diagram

\[\centerline{
\xymatrix{
[\HEE{p,q}{k,l}]^{\ztwo}\,\,\ar[d]^{[\phi]}\ar@{^(->}[rr]^{} & & \HEE{p,q}{k,l}\ar[d]^{\phi} \\
[\hinnfR{p-1,q}{k+1,l}]^{\ztwo}\,\,\ar@{^(->}[rr] & &
\hinnfR{p-1,q}{k+1,l} } }
\]

We have already seen that the map $\phi$ is $\ztwo$-equivariant,
thus it restricts to a continuous map on fixed-points. One can
check with little effort that all parts of the proof of the
quasi-fibration properties still work upon taking fixed points;
that is all homotopies and constructions used are
$\ztwo$-equivariant. The proof therefore applies to the
restriction $[\phi]^{\ztwo}$.

\subsection{Connectivity properties of the spaces}

\label{connectivity}

\begin{prop}
The spaces $\hinnfR{p,q}{k,l}$ are connected for $p,q,k,l\geq 0$
and $k+l\geq 1$. The space $[\hinnfR{p,q}{k,l}]^{\ztwo}$ is
connected for all $k\geq 1$, $l\geq 0$.
\end{prop}

Let $G=\sum \lambda \pi_{G_\lambda}$ be the spectral decomposition of an operator $G\in\hinnfR{p,q}{k,l}$.

Then
\[ H_t(G):= (\frac{1}{1-t})G + (\frac{t}{1-t}) \pi_{\ker (G)} f = \sum_{\lambda\neq 0} (\frac{\lambda}{1-t}) \pi_{G_\lambda} \pm (\frac{t}{1-t}) \pi_{\ker (G)^{\pm}}, \quad t\in [0,1]\]
defines a path from $G$ to the base point using left
multiplication by one generator $f$ of $\ccliff{l,k}$ (we have at
least one such generator by assumption). This splits the
$\ccliff{l,k}-\ccliff{p,q}$-module $\ker (G)=\ker (G^2)$ into two
$\ccliff{p,q}$-modules $\ker (G)^{+} \oplus \ker(G)^{-} = \ker
(G)$ (the eigenspaces of the action of $f$).

Along the path, all eigenvalues move to infinity, symmetrically
about $0$: the eigenvalue $0$ itself splits in two which move out
as well. This move is continuous and depends continuously $t$
(\emph{but not in $G$, see \cite{markert}}).

Note that all properties of $\hinnfR{p,q}{k,l}$ are preserved throughout the path, even though $\ker (G)^{\pm}$ are not $\ccliff{l,k}-\ccliff{p,q}$-modules (their sum is).

\sm

If $G$ is a fixed point of the $\ztwo$-action, i.e. if all
eigenspaces are $\ztwo$-invariant subspaces, we need to use a
generator of the left $\ccliff{l,k}$-action which commutes with
the action $\theta$. This is the case for the generators of
$\ccliff{0,k}$ acting from the left. Our assumptions guarantee
that there is at least one. Then the split parts will again be
$\ztwo$-invariant and all operators in the path will also be fixed
points.

\begin{coro}
The space $\hinnfR{p,q}{k,l}$ is at least $(k+l-1)$-connected. The
space $[\hinnfR{p,q}{k,l}]^{\ztwo}$ is at least $(k-1)$-connected.
\end{coro}

This follows by induction from the deloopings in the previous
section. On fixed-points, the induction stops after $k-1$ steps.

\subsection{Commutation with Bott maps}

\label{commutation}

The diagram
\begin{eqnarray}
\xymatrix{
kr_{n, n}\ar[r]^{\bottoneone}\ar[d]& \Omega^{1,1}kr_{n+1, n+1}\ar[d]\\
KR_{n,n}\ar[r]^{b_{1,1}}& \Omega^{1,1}KR_{n+1, n+1} }
\end{eqnarray}
commutes. This is done by a direct computation.

Notice that the vertical maps in the diagram are essentially
inclusions, forgetting the extra linearity properties of the
operators with respect to  the left Clifford actions. The map
$b_{1,1}$ comes from the structure of $\KR$ as a symmetric
spectrum, as explained in Section \ref{unitmaps}. This map
associates to an operator $F\in \KR_{n,n}$ the map
$$
b_{1,1}(F): \sphere{1,1}  \to \KR_{p+1, q+1};~~ v\mapsto L_v \star
G_{0,0}^{\star 2} \star  F.
$$

On the other hand, the map $\bottoneone$ is given by combining the
formulas (\ref{bott:shift:formula:3}) and
(\ref{bott:shift:formula:4}), i.e. we apply twice the index shift
(formula  (\ref{formula:shift})) and the two maps induced by the
quasi-fibrations. This associates to an operator $F\in \kr_{n,n}$
the map $\bottoneone(F)$:

\begin{eqnarray*}
\sphere{1,1}& \to& \kr_{n+1, n+1}\\  v=(v_1,v_2) &\mapsto&  \widetilde{G_{0,0}}\star\widetilde{G_{0,0}} \star F +  Id_{\dom{\widetilde{G_{0,0}}\star\widetilde{G_{0,0}} \star F}} v_1 \epsilon e_{n+1}  \\
&&+ Id_{\dom{\widetilde{G_{0,0}}\star\widetilde{G_{0,0}} \star F}}
i v_2 \epsilon f_{n+1}
\end{eqnarray*}
which can be rewritten as:

{\small \begin{eqnarray*}
\bottoneone(F)(v) & = & Id_{\ccliff{2,2}}\gradedtensor G_{0,0}^{\star 2}\star F + L_v\gradedtensor Id_{\dom{G_{0,0}^{\star 2}\star F}} \\
 & = & L_v \star G_{0,0}^{\star 2} \star F\\
\end{eqnarray*}
} since multiplication with $v_1 \epsilon e$ and $i v_2 \epsilon
f_{n+1}$ as a right action corresponds to multiplication by $v=
v_1+ i v_2$ as a left action (on the particular copy of
$\ccliff{2,2}$ specified by the generators $e$ and $f$).

\sm This implies that $\kr_{n,n} \to \KR_{n, n}$ models the
$(2n-1)$-connective cover of $\KU_{2n}= \ZZ \times BU$. In fact,
we actually  prove that $\kr_{p, q}\to \KR_{p, q}$ models the
$(p+q-1)$-connective cover of $\ZZ\times BU$ (if $p+q$ is even) or
$U$ (if  $(p+q)$ is odd).

\sm As announced in Section \ref{overview}, we now check that
$\pi_{k}([\covermap _{n,n}]^{\ztwo})$ is an isomorphism for $k\geq
n$. We begin with the observation  \cite{atiyah}, that there is an
equivariant cofibration sequence for all $s,t\geq 0$

$$
\sphere{s,t} \to \sphere{s, t+1}\to\sphere{s+t+1}\wedge \ztwo_+
$$
Applying mapping spaces into a $\ztwo$-space $X$, we get after
taking fixed points the natural fibration sequences of spaces

\begin{eqnarray}
\label{cof:atiy} \Omega ^{s+t+1} \widehat{X} \to  [\Omega ^{s,t+1}
X]^{\ztwo} \to[ \Omega ^{s, t}X]^{\ztwo}
\end{eqnarray}

where $\Omega \widehat{X}$ refers to the non-equivariant loop space of the space $X$ ({\it i.e.} forget all actions).

\sm We note firstly that  $\widehat{\kr}_{p,q}$ is a connective
cover of $\KR_{p, q}$, and that $\KR_{p, q}$ is $\ZZ \times BU$ or
$U$ depending whether $p+q$ is even or not. Hence unequivariantly,
the map $\covermap _{p, q}: \kr_{p, q} \to\KR_{p, q}  $ is
unequivariantly an isomorphism on homotopy groups $\pi_i$ for
$i\geq p+q$. In fact, because our spaces $U$ and $BU$ have
homotopy vanishing every second degree, and more precisely in odd
degrees for $BU$. We therefore even get $\covermap _{n, n}:
\kr_{n, n} \to\KR_{n, n}  $ is unequivariantly an isomorphism on
homotopy groups $\pi_i$ for $i\geq 2n-1$.

\sm Secondly we know that $[\kr_{\bullet, 0}]^{\ztwo} \to
[\KR_{\bullet, 0}]^{\ztwo} $ models the cover map $\ko \to \KO$ so
that  $[\varphi_{n, 0} ]^{\ztwo}$ is an isomorphism on homotopy
groups $\pi_i$ for $i\geq n$. Consider the commutative diagram
$$
\xymatrix{
\Omega^n \widehat{\kr}_{n,n}\ar[r]\ar[d]& [\Omega ^{0,n} \kr_{n,n}]^{\ztwo}\ar[r]\ar[d]& [\Omega ^{0,n-1} \kr_{n,n}]^{\ztwo}\ar[d]\\
\Omega^n \widehat{\KR}_{n,n}\ar[r]& [\Omega ^{0,n}
\KR_{n,n}]^{\ztwo}\ar[r]& [\Omega ^{0,n-1} \KR_{n,n}]^{\ztwo}&. }
$$
The left vertical map is an isomorphism on homotopy groups $\pi_i$
for $i\geq n-1$. The middle one identifies with $\ko_n \to \KO_n $
by the natural equivalences $\Omega^{0,n} \kr_{n,n}\simeq
\kr_{n,0}$ and  $\Omega^{0,n} \KR_{n,n}\simeq \KR_{n,0}$ (plus the
observation that $\KR_{n,0}$ is none but $\KU_0$ acted on by
complex conjugation, and the same for $\kr_{n,0}$, and in this
case fixed point give back the case of $\ko \to \KO$, as in
\cite{markert}), hence is an isomorphism on homotopy groups
$\pi_i$ for $i\geq n$. So we get by the five lemma applied to the
long exact sequence in homotopy of theses fibrations that the
rightmost vertical map is also an isomorphism on $\pi_i$ for
$i\geq n$. Beginning $[\Omega^{0,n}\covermap_{n,n}]^{\ztwo}$, we
induct decreasingly on $k$ to get that $\pi_i
([\Omega^{0,k}\covermap_{n,n}]^{\ztwo})$ is an isomorphism for
$i\geq n$. This gives the desired result for $k=0$.

\end{document}